\magnification=1200
\input amstex
\documentstyle{amsppt}
\NoBlackBoxes
\hsize=6truein
\hcorrection{0.5truein}
\TagsOnRight
\leftheadtext{Boris Apanasov }
\rightheadtext{2-Knots and Representations Of Hyperbolic Groups }
\def\sn{S^n}

\def\hno{\Bbb H^{n+1}}

\define\col{\,:\,}

\def\hra{\hookrightarrow}
\def\lra{\longrightarrow}
\def\ra{\rightarrow}
\def\wh{\widehat}

\def\ov{\overline}
\def\bs{\backslash}
\def\p{\partial}

\def\et{\emptyset}

\define\so{\operatorname{SO}^{\circ}}

\def\hom{\operatorname{Hom}}

\def\rint{\operatorname{int}}
\def\cl{\operatorname{cl}}

\def\isom{\operatorname{Isom}}
\def\mob{\operatorname{M\ddot ob}}
\def\lima{\operatornamewithlimits}

\def\a{\alpha}
\def\da{\delta}
\def\Da{\Delta}

\def\Ga{\Gamma}
\def\La{\Lambda}

\def\th{\theta}
\def\Om{\Omega}

\def\Sa{\Sigma}

\def\th{\theta}

\def\sch{\Cal H}

\def\scq{{\Cal Q}}
\def\scs{{\Cal S}}
\def\scr{\Cal R}
\def\sct{\Cal T}

\define\bh{{\Bbb H}}

\define\N{{\Bbb N}}

\define\br{{\Bbb R}}
\define\R{{\Bbb R}}

\def\fS{\frak S}

\topmatter

\title Two-Dimensional Knots and Representations of Hyperbolic Groups 
\endtitle
\author Boris Apanasov\footnote"\dag"{Supported in part by NSF grant
 DMS-9306311
\hfill \hfill \hfill}
\endauthor
\address Department of Mathematics, University of Oklahoma,
Norman, OK  73019
\endaddress
\email apanasov\@ou.edu \endemail
\abstract\nofrills {ABSTRACT} We describe relations between
hyperbolic geometry and knots of codimension two or, more exactly,
between varieties
of conjugacy classes of discrete faithful representations of the fundamental
groups of hyperbolic $n$-manifolds $M$ into $\operatorname{SO}^{\circ} (n+2,1)$ and
$(n-1)$-dimensional knots in the $(n+1)$-sphere $S^{n+1}$ carrying conservative 
dynamics of hyperbolic lattices $\pi_1(M)$.
This approach allows us to discover a phenomenon of non-connectedness
of these varieties for closed $n$-manifolds $M$, $n\geq 3$, with large enough number of
disjoint totally geodesic surfaces, to construct quasisymmetric infinitely
compounded ``Julia" knots $K\subset S^{n+1}$ which are everywhere wild
and have recurrent $\pi_1(M)$-action, and to
study circle and 2-plane bundles (with geometric 
structures) over closed hyperbolic $n$-manifolds.
\endabstract
\endtopmatter
\bigskip

\document
\head 1. Introduction \endhead
\medskip

This paper studies relations between two spaces which, at first glance,
have nothing in common. The first space is the space of conformal
classes of $m$-knots $K$ of dimension $m\geq 2$ in the $n$-sphere $S^n$, and
the second space is the variety of conjugacy classes of discrete
faithful representations of the fundamental group $\pi_1(M)$ of a
hyperbolic $(m+1)$-manifold $M$ into $SO(n+1,1)$, $n\geq m$, in
particular the Teichm\"uller space $\sct (M)$ of conformal structures on
the manifold $M$.  

Since our knots are at least two-dimensional, those possible relations
are not connected to the problem of existence of hyperbolic metrics
(of constant negative curvature) on knot complements $S^n\bs K$. In fact,
in contrast to 1-knots complements in $S^3$, there are no hyperbolic
metrics on $n$-knots complements in $S^{n+2}$, $n\geq 2$. It follows
from the fact that parabolic ends of
hyperbolic $(n+2)$-manifolds of finite volume are (up to finite cover) products of Euclidean
tori and open intervals, that is hyperbolic metrics may exist only in the
complements of Euclidean surfaces, see \cite{A1} for example.

Nevertheless, one can apply Tukia's \cite{Tu} construction of the
canonical (quasisymmetric) homeomorphisms of the limit sets of isomorphic 
geometrically finite hyperbolic groups to establish such a relation. It
immediately implies (see \cite{A8}): 

{\it There exists a canonical continuous map
$\Phi$ of the Teichm\"uller space $\sct^n(G)$ of geometrically finite
faithful representations of a uniform hyperbolic lattice $G\subset
O(m+1,1)$, $m\geq 2$, to the space $\scq\scs_{m,n}(G)$ of conformal classes of 
$G$-equivariant quasisymmetric embeddings $S^m\hra S^n$. Moreover, this map $\Phi$
is at most a (2-1)-map.}

In particular, if $n=m+2$ the space $\scq\scs_{m,n}(G)$ consists of (conformal 
classes of) $m$-knots $K\subset S^{m+2}$. We also note that the last statement on 
(2-1)-maps follows 
from the Mostow rigidity (\cite{Mo, Ma}) and the Gluck \cite{Gl} theorem (see also
\cite{B, LS, Sw}) that two homeomorphisms of the boundary $\p N(K)\approx S^m
\times S^1$ of a regular neighborhood $N(K)\approx S^{m}\times B^2$
of a $m$-knot $K$, $m\geq 2$, are pseudo-isotopic if and
only if they are homotopic.

Our main results use the map $\Phi$ for a solution of the
connectedness problem for the above spaces of $m$-knots and the
Teichm\"uller spaces associated with a hyperbolic lattice. This solution is based
on our construction in the following theorem.

\proclaim{Theorem 3.1}  For a given nontrivial ribbon $m$-knot
$K\subset S^{m+2}$, $m\geq 2$, there exists a discrete faithful
representation $\rho\col\Ga\ra \mob(m+2)$ of a uniform
hyperbolic lattice $\Ga\subset\isom\bh^{m+1}$ such that the Kleinian group
$G=\rho\Ga$ acts ergodically on the everywhere wild quasisymmetric $m$-knot
$K_\infty=\La(G)\subset\sn$ obtained as an infinite compounding of the
knot $K$.
\endproclaim

As a corollary of this theorem, we have

\proclaim{Theorem 5.2}  Let $\Ga\subset\isom\bh^{m}$ be any uniform hyperbolic
lattice
from the above theorem, and $m\geq 3$.  Then the varieties $\sct^{m+1}(\Ga)$ and 
$\sch^{m+2}(\Ga)$ of conformal and hyperbolic structures on $\Ga$ (respectively, the classes
of discrete faithful representations $\Ga\ra \so (m+2,1)$) are not connected.
\endproclaim

We should mention that this non-connectedness fenomenon has no relation
to non-connectedness of the (smaller) Teichm\"uller space $\sct^{m}(\Ga)$ of
conformal $m$-structures on $M=\bh^m/\Ga$, see
\cite{A5}. Indeed, the topological
obstruction for connectedness of the space $\sct^{3}(M)$ of conformal structures on 
$M=\bh^3/\Ga$  (nontrivial knotting 
of the limit 2-sphere $\La(G)\subset S^3$) is
not an obstruction for connectedness of $\sct^{4}(\Ga)$ and 
$\sch^{5}(\Ga)$.  In fact, the nerve of that
knotting $S^2\hra S^3$ is 1-dimensional, and hence the knotted in $S^3$
2-sphere $\La(G)$ is unknotted in $S^m$, $m\geq 4$.

Another relation between $m$-knots and hyperbolic
manifolds is based on the Gluck's rigidity of Dehn surgeries on trivial
knots of dimension at least 2:

\proclaim{Theorem 6.1}  For a given closed hyperbolic $m$-manifold
$M$, $m\geq 3$, there are at most two non-equivalent circle (or 2-plane) bundles 
over $M$ allowing uniformizable conformal structures or complete
hyperbolic metrics, respectively, with the development on a $(m-1)$-knot complement.
\endproclaim

This paper is preserved in the form it was written in 1995, so it may have no references
to results published after that time. In conclusion, the author would like to thank Scott Carter and 
Masahico Saito for helpful conversations. 

\head 2. Varieties  of representations and geometric structures \endhead
\medskip

Let $M$ be a given closed hyperbolic m-manifold (orbifold),
$m\ge 3$, that is a complete oriented
Riemannian manifold with constant sectional curvature -1, and 
$\pi_1(M)=\Ga\subset \so(m,1)$ its fundamental group isometrically 
acting in the hyperbolic m-space $\bh^m$ as a uniform lattice, $M=\bh^m/\Ga$. 
One can consider the variety
$\scr^n(\Ga)=\hom(\Ga, \so(n,1)$ of all representations of $\Ga$ into 
$\so(n,1), n\ge m$ with the algebraic convergence topology where
the group $\so(n,1)$ acts by conjugations.
Inside of the quotient-variety $\scr^n(\Ga)/\so(n,1)$ there are
two subvarieties

$$\sct^{n-1}(\Ga)\subset\sch^n(\Ga)\subset\scr^n(\Ga)/\so(n,1)\tag2.1$$
which both consist of conjugacy classes of faithful discrete representations,
with an additional condition on representations $\rho$, 
$[\rho]\in\sct^{n-1}(\Ga)$, that 
$\rho(\Ga)\subset \so(n,1)$ are geometrically finite with non-empty 
discontinuity 
sets $\Om(\Ga)\subset S^{n-1}=\p\bh^n$ (the complements of the limit
sets $\La(\Ga)$, $\Om(\Ga)=S^{n-1}\bs \La(\Ga)$). Due to the Sullivan's
stability theorem (see \cite{Su}, \cite{A9, Th.7.2}), the subvariety
$\sct^{n-1}(\Ga)$ is open in $\sch^n(\Ga)$.

There is a natural identification \cite {Mo} 
of the biggest space $\sch^n(\Ga)$ with the space of
$n$-dimensional hyperbolic structures on the group $\Ga$. Such a hyperbolic 
structure on $\Ga$ is
determined (up to inner automorphisms of $\Ga$ and hyperbolic isometries)
by a pair $\{N,\phi\}$, where $N$ is an $n$-dimensional
hyperbolic manifold and $\phi \col \Ga\ra \pi_1(N)$ is an isomorphism. 
In particular, for a closed surface $S_p$ of
genus $p>1$ and  $\Ga=\pi_1(S_p)$, the space $\sch^2(\Ga)$ is the 
Teichm\"uller space $\sct(\Ga)$ of $\Ga$ homeomorphic to $\br^{6p-6}$, 
and the space $\sch^3(\Ga)$ is isomorphic to the product 
$\sct(\Ga)\times\sct(\Ga)$.

Similarly, via holonomy, we have a natural identification  of the space
$\sct^{n-1}(\Ga)$
with the (Teichm\"uller) space of uniformizable $(n-1)$-dimensional conformal
(=conformally flat) structures on $\Ga$ (cf. \cite{A1, A4}). An element of 
this space is determined by a pair
$\{N,\phi\}$, where $N$ is a conformal
(=conformally flat) $(n-1)$-manifold with a non-surjective developing map
$d\col\tilde N\ra S^{n-1}$ of its universal covering space $\tilde N$
into $d(\tilde N)=\Om \subset S^{n-1}$  and $\phi\col\Ga\ra \pi_1(N)$
is a monomorphism corresponding to the short exact sequence

$$0\ra \pi_1(\Om)\ra \pi_1(N)\ra \Ga\ra 0\,.\tag2.2$$

Here, due to \cite{Ka} and \cite{KP}, the developing map $d$ and the 
natural projection $\pi\col\Om\ra \Om/G \cong N$ (where 
$\Ga \cong G=d_{\ast}(\pi_1(N))\subset \so (n,1)$) are
covering maps which factor the universal projection $\tilde N\ra N$.
Two pairs $\{N_0,\phi_0\}$ and $\{N_1,\phi_1\}$ determine the same
conformal structure on $\Ga $ if there is an orientation preserving
conformal homeomorphism $f\col N_0\ra N_1$ such that $f_{\ast}\circ
 \phi_0$
and $\phi_1$ differ (up to the isotropy subgroup $Z(\rho_0)$ of the
inclusion $\rho_0\col\Ga\ra \so (n,1)$) by an inner automorphism of
$\Ga $. (Notice that $f_{\ast}\col\pi_1(N_0)\ra \pi_1(N_1)$ is only well
defined up to inner automorphism of $\pi_1(N_0)$.)

The first surprising results on the varieties $\sct^{n-1}(\Ga)$ and
$\sch^n(\Ga)$ (beyond the Mostow 
rigidity \cite{Mo, Ma}) are that there are large classes of groups $\Ga$
for which these varieties are non-trivial (\cite{A2, A1}), 
but $\sch^n(\Ga)$ is still compact (\cite{T, MS, Mo}). In
particular, due to J.Morgan's weak rigidity theorem \cite{Mo}, we have
(cf. \cite{A4}):

\proclaim{Theorem 2.1} Let $M$ be a closed oriented hyperbolic
m-manifold, $m\ge 3$, with $\pi_1(M)\cong \Ga\subset \so(m,1)$. 
Then, for any $n\ge m$, the
(Teichm\"uller) variety  $\sct^{n-1}(\Ga)$ has a natural
compactification $\overline {\sct^{n-1}(\Ga)}$ such that each of its
points corresponds to a faithful discrete representation
$\rho\col\Ga\ra \so(n,1)$ from the compact variety $\sch^n(\Ga)$.
\endproclaim

\head 3. Knotted $m$-spheres in the $(m+2)$-sphere $S^{m+2}$\endhead
\medskip

An {\it m-knot} $K\subset S^{m+2}$ is an embedding $K\col S^m\hra S^{m+2}$ of
the oriented $m$-sphere into the oriented $(m+2)$-sphere.  Two $m$-knots 
$K_1$ and $K_2$ in $S^{m+2}$ are called {\it equivalent} (or of the same 
{\it type}) if there exists an orientation preserving homeomorphism 
$f\col S^{m+2}\ra S^{m+2}$ such that $f K_1=K_2$.  Obviously, it is an 
equivalence relation,
and we call the equivalence class $[K]$ of a $m$-knot $K\subset S^{m+2}$
the {\it knot type} of $K$.  Those $m$-knots that are equivalent to the natural
inclusion $S^m\subset S^{m+2}$ are called {\it trivial} or {\it unknotted}.

The simplest examples of nontrivial 2-knots in $S^4$ can be obtained by
using the so-called suspensions and spins of classical knots in $S^3$,
see \cite{Ar}, \cite{Ro, \S3J}.  Namely, for a classical 
nontrivial 1-knot $k\subset S^3$, let points $a$ and $b$ lie in disjoint components $S^4_+$
and $S^4_-$ of $S^4\bs S^3$.  Then the join $K$ of the knot $k$ with
$\{a\}\cup\{b\}$ is called {\it the suspension} of the 1-knot $k$, see
Fig. 1.  It holds that $\pi_1(S^4\bs K)\cong\pi_1(S^3\bs k)$, so the
obtained 2-knot $K$ is nontrivial.  

To construct a spun 2-knot
$K\subset S^4$, we consider a classical knot $k$ in the half-space
$\ov{\R^3_+}=\{x\in\R^4\col x_4=0,x_3\geq 0\}$, such that the boundary
2-plane $\R^2=\{x\in\R^4\col x_3=x_4=0\}$ intersects the knot $k$ along a single
arc $\a$.  Then spinning the complementary arc $\beta=\cl (k\bs\a)$
about the plane $\R^2$ sweeps out a 2-sphere $K\subset\R^4$.
Obviously, $K\cap\ov{\R^3}$ consists of a simple loop $k_1$.  Providing an
orientation on $k_1$ which is coherent to that of the arc $\beta\subset k_1$, we
see that $k_1$ is representing the 1-knot that is obtained as the connected sum
$k\#(-k)\subset\R^3$, see Fig. 2.  


The 2-sphere $K\subset S^4$ is oriented so that the orientation of the
2-disk $K\cap S^4_+$ is coherent to that of $k_1$.  So obtained 2-knot
$K\col S^2\hra S^4$ is called the {\it spun 2-knot} of a given 1-knot
$k\subset S^3$, see Fig. 2.  Since every loop in $S^4\bs K$ can be continuously
deformed in $\R^3_+\bs\beta$, we have that

$$\pi_1(S^4\bs K)\cong\pi_1(\R^3_+\bs\beta)\cong\pi_1(\R^3\bs k)\,.$$

We note also that every spun 2-knot $K\subset S^4$ is a 
{\it ribbon 2-knot}. Such 2-knots generalize classical ribbon knots in $S^3$,
see \cite{Ro, BZ}. They can be obtained as follows,
see \cite {Suz} and Fig. 3.


Let $S_0\cup\ldots\cup S_m\subset\R^4$ be a trivial 2-link with $(m+1)$
components (which are trivial non-linked 2-knots) and $f_1,\ldots,f_m$
are appropriate embeddings of 3-balls, $f_i\col [0,1]\times
B^2\hra\R^4$, which make $m$ fusions of the 2-link.  Each of the embeddings
$f_i$ is such that 

$$f_i([0,1]\times B^2)\cap (S_0\cup\ldots\cup S_m)=f_i(\{0,1\}\times B^2)$$
                       has an orientation coherent with that of
the 2-link, and the disks $f_i(\{0\}\times B^2)$ and $f_i(\{1\}\times B^2)$ are
contained in different components of the link.  Then the connected sum
of the spheres $S_0,\ldots,S_m$ and the sphere $f_i(\p([0,1]\times B^2))$
represented by the homological sum

$$(S_0\cup\ldots\cup S_m)+f_i(\p([0,1]\times B^2))=S^1_0\cup\ldots\cup S^1_{m-1}$$
            is a trivial 2-link with $(m-1)$ components.  
Continuing this process of fusions on the link, we finally obtain a 
2-knot which is called a {\it ribbon 2-knot with $m$ fusions}, see Fig. 3.

Similarly, for $m\geq 2$, one can define the above classes of $m$-knots
in $(m+2)$-sphere $S^{m+2}$,
in particular ribbon $m$-knots with a given number of fusions.

Now we can describe a connection between $m$-knots in $S^{m+2}$, 
$m\ge 2$, and varieties of discrete representations of hyperbolic lattices
$\Ga\subset\isom\bh^{m+1}$, which is based on the following theorem.

\proclaim{Theorem 3.1}  For a given nontrivial ribbon $m$-knot
$K\subset S^{m+2}$, $m\geq 2$, there exists a discrete faithful
representation $\rho\col\Ga\ra \mob(m+2)$ of a uniform
hyperbolic lattice $\Ga\subset\so(m+1,1)\cong\isom_+\bh^{m+1}$ such that 
the Kleinian group
$G=\rho\Ga$ acts ergodically on the everywhere wild $m$-knot
$K_\infty=\La(G)\subset\sn$ obtained as an infinite compounding of the
knot $K$.
\endproclaim


\head 4. Basic construction\endhead
\medskip

The proof of the above Theorem 3.1 is based on the author's ``block-building
method'' (see \cite{A3,A5}) and geometrically controlled PL-approximations of
smooth ribbon $(n-2)$-knots $K\subset\sn$, $n\ge 4$.  Namely, we may assume that
the $(n-2)$-dimensional spheres $S_0,\ldots,S_m\subset\sn$ and the embeddings
$f_i\col B^{n-1}\approx[0,1]\times B^{n-2}\hra\sn$ in the definition of a
given ribbon $(n-2)$-knot $K\subset\sn$ are taken in the conformal
category.  That means that all involved spheres are round spheres, and each image
$f_i(B^{n-1})$ is contained in the union of finitely many round 
$(n-1)$-balls $B_j$ in $\sn$, $1\leq j\leq j_i$, such that the boundary spheres of
any two adjacent balls intersect each other along a round 
$(n-3)$-sphere, See Fig 4.


In other words, the $(n-1)$-dimensional ribbon $f_i(B^{n-1})$ can be
obtained from a flat ribbon in $\R^{n-1}$, which is the union of round
balls, by sequential bendings along $(n-2)$-planes. We do that by using
the well known construction of bending deformations (see \cite{A1, A2, Ko, T, A9}).

Let us assume in  addition that, in each round $(n-1)$-ball $B_j$ in the
construction (either a ball from one of the ribbons $f_i(B^{n-1})$ or
one of the balls bounded by spheres $S_k$, $0\leq k\leq m$), there is a discrete 
action of a
hyperbolic group $G_j\subset\isom\bh^{n-1}=\mob(B_j)$.  Up to isotopy of the
$(n-2)$-knot $K$ and the family $\Sa$ of $(n-1)$-balls $B_j$, we may
assume that the groups $G_j$ have bending hyperbolic $(n-2)$-planes
whose boundaries at infinity $\p B_j$ are the intersection spheres $\da_j=\p
B_j\cap \p B_{j+1}$ for the adjacent balls $B_j$
and $B_{j+1}$, and that the stabilizers of $\da_j$ in $G_j$ and
$G_{j+1}$ coincide. We denote such stabilizers by $\Ga_j=G_j\cap G_{j+1}$. 
This property guarantees that the amalgamated free product

$$G=\cdots\underset{\Ga_{j-1}}\to\ast G_j\underset{\Ga_j}\to\ast
G_{j+1} \underset{\Ga_{j+1}}\to\ast\cdots\subset\mob(n)\tag4.1$$
is a Kleinian group isomorphic to a hyperbolic uniform lattice
$\Ga\subset\isom\hno$.

As the result of this geometric construction, we have that our ribbon 
$(n-2)$-knot $K\subset\sn$ is represented as the union $K_0$ of $m+1$
disjoint round $(n-1)$-spheres with $2m$ deleted disjoint round 
$(n-1)$-balls on them,

$$(S_0\cup\ldots\cup S_m)\bs\bigcup^{2m}_{i=1}B^{n-1}_i\,;\quad
 B^{n-1}_i\subset S_k\,,\quad 1\leq i\leq 2m\,,\quad 0\leq k\leq m\,,$$
                        and disjoint 
$(n-2)$-dimensional cylinders corresponding to the ribbons
$f_1,\ldots,f_m$ of the knot $K$.  These disjoint cylinders are the
unions of spherical $(n-2)$-dimensional annuli with disjoint interiors 
which lie on the boundary
spheres $\p B_j$ of the round $(n-1)$-balls $B_j$ in the construction.
Due to the choice of the block-groups
$G_j\subset\mob(B_j)\subset\mob(n)$, the knot $K_0\cong K$ lies in the
interior of the complement of a fundamental polyhedron $F\subset\sn$ of
the product group $G$ in (4.1). Thus $\sn\bs F$ is a regular neighborhood
of $K_0$.  Furthermore, sequentially using bending deformations
(along hypersurfaces corresponding to amalgama-subgroups $\Ga_j$, see \cite{A3, A1}), 
we can construct a hyperbolic uniform lattice
$\Ga\subset\isom\hno$ which has the same amalgamated free product structure as the
group $G$ in (4.1) (and hence, it is isomorphic to $G$) and conformally acts
in one of the $(n-1)$-dimensional balls $B_j$.  We do that in a way
similar to \cite{A3}.

Now a direct application of Tukia's \cite{Tu} isomorphism theorem shows that the
limit set $\La(G)$ of the group $G\cong\Ga\subset\isom\hno$ is an 
$(n-2)$-knot in $\sn$.  We claim that it is the desired everywhere wild knot
$K_\infty\subset\sn$ obtained by the infinite compounding of the knot
$K_0=K$:

$$K_\infty=\ldots\# K\# K \# K \# \ldots\,.\tag4.2$$

Before we go on with the proof, we shall show how this construction
works in the case of 2-knots.  We illustrate it by a simplest ribbon 
2-knot $K\subset \R^4$ obtained by one fusion from two unlinked 2-spheres
$S_0$ and $S_1$, see Fig. 3.  This knot $K$ is also the spun 2-knot of
the classical trefoil knot $k\subset\R^3$, see Fig. 2.

We can take the knot $K\subset S^4$ to be a PL-knot shown in Fig. 5, that
is as a ribbon 2-knot obtained by one fusion of two unlinked boundaries
of 3-dimensional cubes $Q_0$ and $Q_1$ in 3-planes $\R^3\times\{0\},\,
\R^3\times\{54\}\subset \R^4$, respectively.

To have better geometric control on the blocks that we are going to use in the
construction, we shall use (instead of the spheres $S_0$ and $S_1$ in the 
definition of a ribbon knot)
the cubes $\p Q_i$ with edges parallel to the
coordinate axes in $\R^4$.  As the ribbon
$f_1\col B^3\hra\R^4$ we shall also use the union of 3-dimensional
(smaller) cubes $Q_j$, $2\leq j\leq m$, in $\R^4$ as it is indicated in
Fig. 5.


Here the edges of the equal small cubes $Q_j$, $j\geq 2$, are parallel to the
coordinate axes in $\R^4$. Their size shall be later determined by the size of 
the first two bigger cubes
$Q_0\cong Q_1$. The plus signs in Fig. 5
show the character of intersections in our 3-dimensional projection of
$K$ of the boundaries $\p Q_0$ and $\p Q_1$ with the boundary of 
3-dimensional tube which is the union $\bigcup_{2\leq j\leq m}Q_j$ of
small cubes.  Later we shall also explain our choice of parallel 
3-planes $\R^3\times\{-27\}$, $\R^3\times\{0\}$, $\R^3\times\{54\}$ and
$\R^3\times\{81\}$ in $\R^4$ which contain some of the cubes $Q_j$ and are
orthogonally joined by tubes that are boundaries of the union of the
remaining small cubes.

Now we define discrete block-groups $G_j$ associated with the cubes $Q_j$.
Although $G_j$ are isomorphic to hyperbolic isometry groups in $\bh^3$, it is 
more convenient to use quasi-Fuchsian deformations of these hyperbolic groups 
(bendings as in cite{A2, A3}) so
that the obtained groups $G_j$ match the cubes $Q_j$ in the following sense.

Assuming $K=K_0$, we cover the 2-knot 
$K\subset\bigcup_{0\leq j\leq m}\p Q_j$ by a family
$\Sa=\{b_{ji}\}$ of closed round 4-balls $b_{ji}$ whose boundary spheres $\p
b_{ji}$ are orthogonal to $K$.  Namely, in the first step, we take 
4-balls $b_{ji}$ centered at the vertices of the cubes $Q_j$, $0\leq j\leq m$,
whose radii $r_{ji}$ are equal to each other if either $j=0,1$ or $2\leq
j\leq m$.  One more condition on these radii $r_{ji}$ is that
$b_{ji}\cap b_{kl}\neq \et$ only if the centers of the different balls
$b_{ji}$ and $b_{kl}$ are the ends of a common 1-edge of one of the
cubes.  In the latter case, the magnitude of the exterior dihedral angle
bounded by the spheres $\p b_{ji}$ and $\p b_{kl}$ should equal $\pi/3$.
These 4-balls $b_{ji}$ do not cover the entire knot $K$.  On each square
2-side $X\subset\p Q_j\cap K$, we have uncovered 4-gon bounded by
circular arcs.  We cover such a 4-gon by five additional 4-balls
$b_{ji}$ centered at $X$ and whose boundary spheres $\p b_{ji}$
intersect (orthogonally) only those previously constructed balls that are
centered at the vertices of $X$.  Among these five new balls, the first
four sequentially intersect each other with the external 
dihedral angles $\pi/3$. The fifth ball is centered at the center of $X$ 
and (orthogonally) intersects 
only the last new four balls, see Fig. 6.  After that, we still have
uncovered those two 2-sides $X_0$ and $X_1$ of the big cubes $Q_0$ and
$Q_1$ that are (orthogonally) joined by the tube 
$\bigcup_{2\leq j\leq m}Q_j$.  Here we assume that $Q_0\cap Q_2$ and 
$Q_1\cap Q_m$ are small
squares centered at the centers of $X_0$ and $X_1$, respectively.
Furthermore, we choose {\it the size} of the cubes $Q_j$ so that $Q_j$, $j\geq 2$, 
are unit cubes and the cubes $Q_0$ and $Q_1$ have
the size which matches the covering family $\{b_{ji}\}$.  In fact, the
boundary spheres of the four additional balls centered at
$\rint (X_0)$ orthogonally intersect the corresponding spheres
centered at the four vertices of the small cube $Q_2$, see Fig. 6.


This completes the construction of the family
$\Sa=\{b_{ji}\}$ of 4-balls covering the knot $K$.  The union of these
balls, $\bigcup_{i,j}\text{int}b_{ji}=N(K)$, is a regular neighborhood
of the PL-ribbon knot $K$.

Notice that we can take the size of cubes $Q_j$, $j\geq 2$,
arbitrarily smaller than that of the cubes $Q_0$ and $Q_1$. To do that,
we repeat the above process of covering the sides $X_0$ and $X_1$ by
balls $b_{ji}$ where, instead of the vertices of $X_0$ (and $X_1$),
we take the centers of the four new small balls. Then each of the annuli
in $X_0\bs Q_2$ and $X_1\bs Q_m$ will be covered by
$(4+8k)$ additional balls $b_{ji}$ (for sufficiently large integer
$k\geq 0$) instead of the above four additional balls corresponding to
$k=0$.  This allows us to take the ribbon $f_1\col B^3\hra\R^4$ as thin
as we need.  

We define a discrete block-group $G_j$ associated with a cube $Q_j$,
$0\leq j\leq m$, as the group generated by reflections with respect to
all spheres $\p b_{ji}$, that is, with respect to all spheres
$\p b_{ji}$ that intersect the cube $Q_j$.  Obviously, $G_j$ is discrete
because all spherical dihedral angles with edges $\p b_{ji}\cap \p b_{jl}$ are either
$\pi/3$ or $\pi/2$ (see \cite{A1, Mas} for example).  Furthermore,
$G_j$ preserves each of (coordinate) 3-planes $\R^3\subset\R^4$ that contain the cube
$Q_j$.  In such a 3-plane $\R^3$, the group $G_j$ can be deformed by
bendings to a Fuchsian group acting in a 3-ball $B^3\subset\R^3$.  That
is why we can consider the groups $G_j$ as discrete subgroups in
$\isom\bh^3$, $G_j\cong G'_j\subset\isom\bh^3$.

For any two adjacent cubes $Q_j$ and $Q_{j+1}$, the groups $G_j$ and
$G_{j+1}$ have a common subgroup $\Ga_j=G_j\cap G_{j+1}$ which is
generated by four reflections with respect to the spheres centered at
the vertices of the square $Q_j\cap Q_{j+1}$.  So we can apply the
Maskit combination \cite{Mas, A1} and obtain a Kleinian group
$G\subset\mob(4)$ as the free amalgamated product in (4.1).  For the
group $G$, we can take the complement of a regular neighborhood $N(K)$
of the knot $K$ to be a fundamental polyhedron $P\subset S^4$:

$$P=\ov{\R^4}\bs N(K)\,,\quad N(K)=\bigcup_{i,j}\rint b_{ji}\,.\tag4.3$$

We remark that, for each amalgamated free product
$G_j\underset{\Ga_j}\to\ast G_{j+1}$, we can use a bending deformation
along the hyperbolic 2-plane $H_j$ whose boundary circle $\p H_j$ is the
limit circle of the amalgama subgroup $\Ga_j$.  As a result, we get a
new hyperbolic isometry group $G'_j\subset\isom\bh^3$ which is isomorphic
to $G_j\underset{\Ga_j}\to\ast G_{j+1}$.  Applying this process $m$
times, we obtain a cocompact discrete group $\Ga\subset\isom\bh^3$
isomorphic to the group $G$.  

In dimension $n=4$, there is another (non-algorithmical) way to get such a 
unique hyperbolic lattice $\Ga$ by
using the Andreev-Rivin classification of hyperbolic compact polyhedra
in $\bh^3$.  Namely, the boundary $\p P$ of the polyhedron in (4.3) has
the combinatorial type of $S^2\times S^1$ where the 2-sphere $S^2$ is
decomposed into the union of spherical polygons.  In fact, $\p P$ is the
union of 3-sides each of which is the annulus on a sphere $\p b_{ji}$,
i.e. each 3-side is the product of a spherical 2-polygon $D_{ji}$ and
the circle $S^1$.  The dihedral angles between such 3-sides are determined by the
corresponding 3-dimensional dihedral angles bounded by 2-spheres  $\p
b_{ji}\cap\R^3$ in the corresponding 3-planes $\R^3\subset\R^4$, so they
are either $\pi/3$ or $\pi/2$, and the Andreev-Rivin conditions apply
(see \cite{An, Ri}).  It follows that the combinatorial type of the
4-polyhedron $P$ determines the combinatorial type of a 3-dimensional
compact hyperbolic polyhedron $P'\subset\bh^3$, with the same magnitudes
of dihedral angles as those for $P$.  Thus the group $\Ga\subset\isom\bh^3$
generated by reflections in sides of $P'$ is a uniform hyperbolic
lattice isomorphic to the group $G$.

\remark{Remark 4.1}  The above observation that we can take the ribbon
$f_1\col B^3\hra\R^4$ ``arbitrarily thin'' makes it possible to apply
the above block-groups $G_j\subset\mob(4)$, $G_i\cong
G'_j\subset\isom\bh^3$, to represent an arbitrary ribbon 2-knot 
$K\subset S^4$ as the knot which lie on the boundary of the union of 
3-cubes similar the above cubes $Q_j$. 
\endremark

To finish the proof of Theorem 3.1, we need to show that the knot 
$K_\infty=\La(G)$ in (4.2) is an everywhere wild $(n-2)$-knot if the knot $K$ 
is nontrivial. 

Let $\fS=G(\bigcup_{ij}\p b_{ji})$ be the $G$-orbit of the boundary
spheres $\p b_{ji}$ of the balls in the covering $\Sa$ of the PL-knot
$K$.  We can use the word norm $|g|$ of elements $g\in G$ with respect
to generators of $G$, which are reflections in sides of the fundamental
polyhedron $P=\sn\bs N(K)$ in (4.3), to define a partial ordering on
$\fS$.  Namely, for two spheres $S_1,S_2\in\fS$, we say $S_2\succcurlyeq S_1$ 
if $\rint S_2\subseteq\rint S_1$.  It allows us to enumerate the
set $\fS$ by a bijection $q\col\N\ra\fS$ so that it is compatible with
the ordering of $\fS$, that is the map $q^{-1}$ preserves this partial
order.  Then we have a nested sequence of compacta,

$$P_0=P\subset P_1=P_0\cup g_1(P_0)\subset\ldots\subset P_k=P_{k-1}\cup
g_k(P_{k-1})\subset\ldots\subset\Om(G)\,,\tag4.4$$
                        where the elements $g_i\in G$ are the reflections 
with respect to $i$-th
spheres $S_i\in\fS$ each of which contains a side of the $(i-1)$-th
polyhedron $P_{i-1}$.  

The complement $\sn\bs P_i$ of each of the compacta $P_i$ in (4.4) is a
regular open neighborhood of an $(n-2)$-knot $K_i$ which is obtained from
the knot $K$ by sequential connected sums:

$$K_0=K,\quad K_1=K_0\# K_0,\ldots,K_i=K_{i-1}\# K_{i-1},
\ldots\,.\tag4.5$$

Since the limit set $\La(G)=\sn\bs\Om(G)$ is homeomorphic to the limit set
$\La(\Ga)=S^{n-2}$ (Tukia's \cite{Tu} isomorphism theorem),  
$\La(G)$ is an embedded $(n-2)$-sphere in $\sn$.  
We denote $\La(G)=K_\infty$ and claim that it is an everywhere wild 
$(n-2)$-knot in $\sn$.  Obviously, $K_\infty=\bigcap_i\ov{N(K_i)}$
where, for any $i$, 
$N(K_i)=\sn\bs P_i$ is a  regular neighborhood of the knot $K_i$
in (4.5). Due to (4.4), the nested sequence $\{\ov{N(K_i)}\}$ is
decreasing to its intersection, $K_\infty$.

Due to the Alexander duality \cite{Sp} applied to
$\Om(G)=\sn\bs K_\infty$, we have that 
$H_1(\Om(G);\Bbb Z)\cong H^{n-2}(S^{n-2};\Bbb Z)=\Bbb Z$.  
Thus we can consider an infinite cyclic
covering space $\wh\Om$ of $\Om(G)$.  Now we are concerned with 
the integral
homology $H_\ast(\wh\Om)=H_\ast(\wh\Om;\Bbb Z)$ with $\La$-module
structure where $\La$ denotes the ring of finite Laurent polynomials
with integer coefficients.  Namely, choosing a generator
$\tau\col\wh\Om\ra\wh\Om$ of the deck transformation group of the cyclic
covering $\wh\Om\ra\Om(G)$, we define the product of an element
$p(t)=\lima{\sum}_{-r\leq i\leq s}c_it^i\in\La$ with an element 
$\a\in H_j(\wh\Om)$ as

$$p(t)\cdot\a=\sum^s_{i=-r}c_i\tau^i_\ast\cdot\a\in H_j(\wh\Om)\,.$$
Here $\tau_\ast\col H_j(\wh\Om)\ra H_j(\wh\Om)$ is the homology
isomorphism induced by $\tau$.  Thus it defines the $\La$-module
$H_\ast(\wh\Om)$ which is known as the Alexander invariant of the knot
$K_\infty\subset\sn$.  As a shorthand description of the first homology
$H_1(\wh\Om)$ of the infinite cyclic covering of the knot $K_\infty$
complement $\Om(G)$, one can also use its Alexander polynomial
$\Da_{K_\infty}(t)$, see \cite{Ro, Ch. 7}.

\proclaim{Lemma 4.2}  Let $\wh\Om$ and $\wh P_k$, $k\geq 0$, be infinite
cyclic coverings of $\Om(G)$ and $P_k$ in (4.4).  Then we have a nested
sequence 
$$\wh P_0\subset\wh P_1\subset\ldots\subset\wh
P_k\subset\ldots\wh\Om\,.$$
\endproclaim

\demo{Proof}  The nested sequence in (4.4) defines a sequence of
monomorphisms of the fundamental groups as follows:

$$\pi_1(P_0)\hra\pi_1(P_1)\hra\ldots\hra\pi_1(P_k)\hra
\pi_1(P_{k+1})\hra\ldots\hra\pi_1(\Om(G))\,.$$

%
%
%
%
\catcode`\@=11 
\def\hookrightarrowfill{$\m@th\mathord\lhook\mkern-3mu\mathord-\mkern-6mu%
  \cleaders\hbox{$\mkern-2mu\mathord-\mkern-2mu$}\hfill
  \mkern-6mu\mathord\rightarrow$}
\def\hookleftarrowfill{$\m@th\mathord\leftarrow\mkern-6mu%
  \cleaders\hbox{$\mkern-2mu\mathord-\mkern-2mu$}\hfill
  \mkern-6mu\mathord-\mkern-3mu\mathord\rhook$}
\atdef@ C#1C#2C{\ampersand@\setbox\z@\hbox{$\ssize
 \;{#1}\;\;$}\setbox\@ne\hbox{$\ssize\;{#2}\;\;$}\setbox\tw@
 \hbox{$#2$}\ifCD@
 \global\bigaw@\minCDaw@\else\global\bigaw@\minaw@\fi
 \ifdim\wd\z@>\bigaw@\global\bigaw@\wd\z@\fi
 \ifdim\wd\@ne>\bigaw@\global\bigaw@\wd\@ne\fi
 \ifCD@\hskip.5em\fi
 \ifdim\wd\tw@>\z@
 \mathrel{\mathop{\hbox to\bigaw@{\hookrightarrowfill}}\limits^{#1}_{#2}}\else
 \mathrel{\mathop{\hbox to\bigaw@{\hookrightarrowfill}}\limits^{#1}}\fi
 \ifCD@\hskip.5em\fi\ampersand@}
\atdef@ D#1D#2D{\ampersand@\setbox\z@\hbox{$\ssize
 \;\;{#1}\;$}\setbox\@ne\hbox{$\ssize\;\;{#2}\;$}\setbox\tw@
 \hbox{$#2$}\ifCD@
 \global\bigaw@\minCDaw@\else\global\bigaw@\minaw@\fi
 \ifdim\wd\z@>\bigaw@\global\bigaw@\wd\z@\fi
 \ifdim\wd\@ne>\bigaw@\global\bigaw@\wd\@ne\fi
 \ifCD@\hskip.5em\fi
 \ifdim\wd\tw@>\z@
 \mathrel{\mathop{\hbox to\bigaw@{\hookleftarrowfill}}\limits^{#1}_{#2}}\else
 \mathrel{\mathop{\hbox to\bigaw@{\hookleftarrowfill}}\limits^{#1}}\fi
 \ifCD@\hskip.5em\fi\ampersand@}
\catcode`\@=13
Furthermore, we have the following commutative diagram:

$$\CD
\cdots \pi_1(P_{k-1}) @CC{i_k}C \pi_1(P_k) @CCC \cdots @CC i C
\pi_1(\Om(G))\\
 @VVV  @VVV @. @VVV \\
\cdots H_1(P_{k-1}) @C{i_{k\ast}}CC H_1(P_k) @CCC \cdots @C{i_\ast}CC
H_1(\Om(G))\\
 @|  @| @. @| \\
\cdots\phantom{X} \Bbb Z\phantom{P_{k-1})} @>\cong>> \Bbb Z @>\cong>> \cdots @>\cong>> \Bbb
Z\endCD$$
Here the vertical maps correspond to the Abelinization.
So the lemma follows.
\enddemo

\proclaim{Lemma 4.3}  $H_1(\wh\Om;\Bbb Z)\neq 0$.
\endproclaim

\demo{Proof}  Due to Lemma 4.2, we have excisive triads defined by $X=\wh
P_i$ and $Y=\wh\Om\bs\wh P_k$.  So we have a Mayer-Vietoris exact
sequence as follows:
$$\CD
\cdots @>>> H_2(X\cap Y) @>>> H_2(X)\oplus H_2(Y) @>>> H_2(X\cup Y)
@>>>\\
\hphantom{\cdots }@>>> H_1(X\cap Y) @>>> H_1(X)\oplus H_1(Y) @>>>
H_1(X\cup Y) @>>> \cdots\endCD$$
where $X\cap Y=\p\wh P_k$ and $X\cup Y=\wh\Om$.

For each $k\geq 0$, the boundary $\p P_k$ is homeomorphic to the product
$S^{n-2}\times S^1$.  It implies that $\p\wh P_k\approx S^{n-
2}\times\R$, and hence $H_1(\p\wh P_k)=0$.  Therefore, the map
$$H_1(\wh P_k)\oplus H_1(\wh\Om\bs\wh P_k)\lra H_1(\wh\Om)$$
is injective.

On the other hand, $H_1(\wh P_k)$ is not trivial due to the initial
condition that the $(n-2)$-knot $K\subset\sn$ is nontrivial.  This
completes the proof.
\enddemo

To finish the proof of Theorem 3.1, we observe that the $(n-2)$-knot
$K_\infty=\La(G)\subset\sn$ is nontrivial due to nontriviality of its
Alexander invariant $H_\ast(\wh\Om)$, which follows from Lemma 4.3.  On the
other hand, $K_\infty$ is invariant for a non-elementary Kleinian group
$G\subset\mob(n)$, and hence it is a wild knot due to \cite{Ku}.
The latter fact also follows from Lemma 4.3 and the additivity of the
Alexander invariant $H_\ast(\wh\Om)$ with respect to connected sum 
(4.5) of knots, see
\cite{L}.  Obviously, any point $z\in K_\infty=\La(G)$ that is
the attractive fixed point of a loxodromic element $g\in G$ is a wild
point of the knot $K_\infty$.  The proof is completed by the well known
fact (see \cite{A1} for example) that such loxodromic fixed points are 
dense in the limit set $\La(G)=K_\infty$.

\head 5. Quasisymmetric Embeddings and Variety Components\endhead
\medskip

Let $\scq\scs_{m,n}$,  $ n\geq m\geq 2$, be the space  of conformal classes 
of quasisymmetric 
embeddings $S^m\hra S^n$ with the compact-open topology. Up to M\"obius 
transformations, such 
quasisymmetric embeddings $f\col S^m\hra S^n$ of $m$-sphere 
$S^m=\br^m\cup{\infty}$ can be represented by 
quasisymmetric embeddings of Euclidean $m$-space, 
$ f\col \br^m\hookrightarrow \br^n$, which satisfy the following
condition \cite{TV}:

$$\frac {1}{\eta(\rho)}\leq \frac {|f(y)-f(x)|}{|f(z)-f(x)|}\leq
\eta(\rho)\quad \text{if} \quad
\frac {1}{\rho}\leq \frac {|y-x|}{|z-x|}\leq \rho,\,\, \rho>0\,,$$
where $\eta\col [0,\infty)\ra [0,\infty)$ some homeomorphism and 
$x,y,z\in\br^m$.

For a given lattice $\Ga\subset O(m+1,1)$, we have a subspace
$\scq\scs_{m,n}(G)\subset \scq\scs_{m,n}$ consisting of all
$\Ga$-equivariant quasisymmetric embeddings $f\col S^m\hra S^n$.
Here the groups $\Ga$ and  $f\Ga f^{-1}$ act on 
$S^m=\p H^{m+1}$ and $f(S^m)\subset S^n$ by restrictions of 
M\"obius transformations in the corresponding spheres. 

Since the canonical homeomorphisms of the limit sets in Tukia's \cite{Tu} 
isomorphism theorem are in fact quasisymmetric, we immediately have an 
additional metric property of the knot $K_\infty$ in Theorem 3.1:

\proclaim{Corollary 5.1}  For a given nontrivial ribbon $m$-knot
$K\subset S^{m+2}$, $m\geq 2$, there is a quasisymmetric embedding $f\col
S^{m}\hra S^{m+2}$ whose image is an everywhere wild $m$-knot
$K_\infty=f(S^{m+2})$, infinitely compounded from $K$.
\endproclaim

Considering conjugations of the inclusion of a lattice $\Ga\subset \so(m+1,1)$
to $\so(n+1,1)$ by quasiconformal self-homeomorphisms
of $\bh ^{n+1}$ compatible with the $\Ga$-action, we have quasiconformal
deformations of the inclusion $\Ga\subset \so(n+1,1)$, i.e. curves in 
the varieties of representations $\sct^{n}(\Ga)\subset\sch^{n+1}(\Ga)$
corresponding
to continuous families of such conjugations. Following to the classical
terminology, we call such representations as quasi-Fuchsian ones.
The Sullivan's stability theorem
\cite{Su, A9} implies that the set of conjugacy classes of quasi-Fuchsian
representations is an open connected subspace of $\sct^{n}(\Ga)\subset
\sch^{n+1}(\Ga)$. We denote the connected components of $\sct^{n}(\Ga)$
and $\sch^{n+1}(\Ga)$ containing this open subspace by $\sct^{n}_{\circ}(\Ga)$ and 
$\sch^{n+1}_{\circ}(\Ga)$. 

Surprisingly, in
contrast to the classical Teichm\"uller theory and the trivial space
$\sch^{m+1}(\Ga)$, for $m\ge 2$,  the varieties $\sct^{3}(\Ga)$ and 
$\sch^{4}(\Ga)$ may be disconnected for hyperbolic 3-lattices $\Ga\subset\so(3,1)$. 
Indeed, as it was pointed out by the author \cite{A5},
for some faithful discrete representations $\rho\col \Ga\ra \so(4,1)$, 
the limit set $\La(\rho(\Ga))$ may be an everywhere wild
2-sphere in $S^3$ (the boundary sphere
of a wildly embedded 3-ball). On the other hand, this is not a reason
for $\sct^{4}(\Ga)$ and $\sch^{5}(\Ga)$  to be disconnected because 
the nerve of the constructed in \cite{ A5} knotting is 1-dimensional, so
the limit 2-sphere $\La(\rho\Ga)\subset S^3$ is unknotted in $S^4$. 

Nevertheless, as another application of Theorem 3.1, we have:

\proclaim{Theorem 5.2} Let $\Ga\subset \so(k,1)$ be any uniform lattice
from Theorem 3.1, and $k\ge 3$. Then the varieties $\sct^{k+1}(\Ga)$ and 
$\sch^{k+2}(\Ga)$ of conformal and hyperbolic structures on $\Ga$
(or equivalently, of conjugacy classes of discrete faithful representations 
$\Ga\ra \so (k+2,1)$)
are disconnected.
  \endproclaim

\demo{Proof} It immediately follows from Theorem 3.1 and the Sullivan's
stability theorem because both (isomorphic) groups 
$\Ga\subset \isom \bh^k\subset\isom \bh^{k+2}$ and
$G=\rho\Ga\subset \isom \bh^{k+2}$ constructed in Section 4 are convex
cocompact (geometrically finite loxodromic groups). Namely,
we see that the
''knotted" representation $\rho$ does not belong to the closure of 
the quasi-Fuchsian
components $\sct^{k+1}_{\circ}(\Ga)$ and $\sch^{k+2}_{\circ}(\Ga)$ 
because, for any quasi-Fuchsian representation $\rho'$, the limit set 
$\La(\rho'\Ga)$ is an unknotted 
$(k-1)$-sphere in $S^{k+1}$. 
\enddemo

  \head 6. Fiber Bundles Over Hyperbolic Manifolds \endhead
\medskip

We conclude our paper by pointing out a significant 
difference between geometric structures on fiber bundles
over hyperbolic surfaces and $m$-manifolds, $m\ge 3$, which is based on
Gluck's rigidity for Dehn
surgery on high-dimensional knots (in contrast to Dehn surgery on
classical 1-knots in $S^3$). 

Given a closed hyperbolic $m$-manifold $M$ with
$\pi_1(M)=\Ga\subset \so(m,1)$, each faithful discrete representation
$\rho\col\Ga\ra\so(m+2)$ (a conjugacy class $[\rho]\in \sch^{m+2}$)
defines a $(m+2)$-dimensional hyperbolic manifold $E=\bh^{m+2}/\rho\Ga$ 
and a closed conformal $(m+1)$-manifold $N$ at infinity of $E$,
$N=\Om(\rho\Ga)/\rho\Ga$, whose fundamental group satisfies 
the exact sequence (2.2) with $G=\rho\Ga\subset \so (m+2,1)$.
The manifolds $E$ and $N$
themselves are 2-plane and $S^1$ bundles over $M$
provided the exterior of the $(m-1)$-knot $\La(\rho\Ga)\subset S^{m+1}$
is homeomorphic to $B^m\times S^1$, that is the knot
is trivial \cite{ Go}).  
Moreover, it turns out that (besides trivial bundles 
corresponding to quasi-Fuchsian representation of $\Ga=\pi_1(M)$)
there are such non-trivial circle and 2-plane bundles over $M$ 
having conformal and hyperbolic
structures and realized by the manifolds $E$ and $N$, correspondingly. 
Namely, Gromov, Lawson and Thurston \cite{GLT} have produced startling 
constructions of such representations $\rho$ whose hyperbolic and conformal
manifolds $E$ and $N$ are non-trivial 2-plane and circle bundles over a 
closed hyperbolic surface, see also \cite{Ku}. Those fiber bundles are related to conformal 
realizations of Dehn surgeries on a classical 1-knot $S^1\hra S^3$.
Similarly, existence of a hyperbolic structure on a 2-plane bundle over a hyperbolic 
$m$-manifold $M$ (conformal structure on a circle bundle over $M$) whose
holonomy group $\rho\Ga$, $\Ga=\pi_1(M)$, has an $(m-1)$-knot in
$S^{m+1}=\p\bh^{m+2}$ as the limit set implies a conformal realization
of the corresponding Dehn surgery on the trivial knot 
$S^{m-1}\subset S^{m+1}$. 

 However, the Dehn surgery on
high-dimensional knots $S^{m-1}\hra S^{m+1}$, $m\ge 3$ is very rigid. 
This is related to the fact, for the first time observed by \cite{Gl} in dimension
$m=3$ (see also \cite{B, LS, Sw}), that two homeomorphisms of the 
boundary $\p N(K)\approx S^{m-1}\times S^1$ of a regular neighborhood 
$N(K)\approx S^{m-1}\times B^2$
of a $(m-1)$-knot $K\subset S^{m+1}$, $m\geq 3$, are pseudo-isotopic if and
only if they are homotopic.  The group of pseudo-isotopy classes of
homeomorphisms of $S^{m-1}\times S^1$ is thus isomorphic to $\Bbb
Z_2\times\Bbb Z_2\times\Bbb Z_2$.  Here the first two factors correspond
to orientation-reversals of $S^{m-1}$ and $S^1$ respectively, and the
third is generated by the following homeomorphism $\eta\col S^{m-1}\times
S^1\ra S^{m-1}\times S^1$,

$$\eta(x,\th)=(\tau(\th)(x),\,\th)\,;\quad x\in S^{m-1}\,,\quad \th\in
S^1\,,\tag6.1$$
              where $\tau(\th)$ is the rotation of the sphere $S^{m-1}$ about its
polar $S^{m-3}$ through the angle $\th$.

Therefore, in contrast to the classical 1-knots in $S^3$, each 
$(m-1)$-knot $K\subset S^{m+1}$, $m\geq 3$, has the only one nontrivial 
Dehn surgery
which is determined by the homeomorphism (6.1).  This makes fiber bundles
over closed hyperbolic $m$-manifolds, whose fibers are either 
2-planes or circles and which have either hyperbolic
structures (of infinite volume) or conformal structures,  correspondingly,
 more rigid than the analogous fibrations over hyperbolic surfaces, 
see \cite{GLT}.  In particular it implies \cite{A8}:

\proclaim{Theorem 6.1}  For a given closed hyperbolic $m$-manifold
$M$, $m\geq 3$, there are at most two non-equivalent circle (or 2-plane) bundles 
over $M$ allowing uniformizable conformal structures or complete
hyperbolic metrics, respectively, with the development on a $(m-1)$-knot complement.
\endproclaim

Finally we remark that there are at most finite number of equivalence
classes of conformal $(m+1)$-manifolds $N$ (hyperbolic $(m+2)$-manifolds 
$E,\ \p E=N$) homotopy equivalent to a given closed hyperbolic $m$-manifold
$M$ whose developments are onto a (nontrivial) $(m-1)$-knot complement.
The number of such manifold equivalence classes depends of geometry
of $M$, especially, of the number $b(M)$ of disjoint totally geodesic
surfaces in $M$. In our example in Section 4 related to the (trefoil) 2-knot $K\subset S^4$, 
this number is more than 100 and, we think, it cannot be less than 
$C_m\cdot e$ where $C_m>1$ is an universal constant and $e=e(K)$ is the 
potential energy of the m-knot $K$ \cite{ BFHW, AS} (actually, for the used trefoil $k\subset S^3,\
 e(k)\approx 74$).

\vfil
\newpage
\def\ref#1{[#1]}
\eightpoint
\parindent=36pt

\head  REFERENCES\endhead
\bigskip

\frenchspacing

\item{\ref{An}} E.M. Andreev, On convex polyhedra in Lobachevskii
space. - Math. USSR, Sbornik {\bf 10} (1970), 413-440.

\item{\ref{A1}}  Boris Apanasov,   Discrete  groups   in  Space   and
Uniformization  Problems.- Math. and  Appl. {\bf 40},  Kluwer
Academic Publishers, Dordrecht, 1991, IX + 484 pp.

\item{\ref{A2}} \underbar{\phantom{Apanas}}, Nontriviality  of Teichm\"uller space for
Kleinian  group in  space. - Riemann Surfaces  and Related
Topics: Proceedings  of the 1978 Stony  Brook Conference, I.Kra  
and B.Maskit, Eds  (Ann. of  Math. Studies {\bf 97}) -
Princeton Univ. Press, 1981, 21-31.

\item{\ref{A3}}  \underbar{\phantom{Apanas}}, Quasisymmetric embeddings of
a closed ball inextensible in neighborhoods of any boundary point. - 
Ann. Acad. Sci. Fenn. Ser. A I Math. {\bf 14} (1989), 243-255.

\item{\ref{A4}}  \underbar{\phantom{Apanas}}, Hyperbolic cobordisms  and
conformal structures.- `` Topology'90", W. de
Gruyter-Verlag, Berlin-New York, 1992, 27-36.

\item{\ref{A5}}  \underbar{\phantom{Apanas}},  Non-standard uniformized
conformal structures on hyperbolic manifolds.- 
  Invent. Math. {\bf 105} (1991), 137-152.

\item{\ref{A6}} \underbar{\phantom{Apanas}},  Deformations of conformal structures on
hyperbolic manifolds. - J. of Diff. Geom. {\bf 35}
(1992), 1-20.

\item{\ref{A7}} \underbar{\phantom{Apanas}}, Varieties of discrete
representations of hyperbolic lattices.- Algebra (Yu.L.Ershov a.o., Eds),
M.I.Kargapolov Memorial Volume, W. de Gruyter-Verlag, 
Berlin-New York, 1996, 7-20.

\item{\ref{A8}} \underbar{\phantom{Apanas}},  Quasi-symmetric knots and 
Teichm\"uller spaces. -  Russian Acad. of Sci. Dokl. Math., {\bf 55:3}
(1997), 319-321.

\item{\ref{A9}} \underbar{\phantom{Apanas}}, Conformal Geometry of 
Discrete Groups and Manifolds. - De Gruyter Expositions in Math. {\bf 32},
W. de Gruyter, Berlin-New York, 2000, XIV+523 pp.

\item{\ref{Ar}} E.~Artin, Zur Isotopie zweidimensionaler Fl\"achen 
im $\R^4$. - Hamburg Abh. {\bf 4} (1925), 174--177.

\item{\ref{AS}} David Auckly and Lorenzo Sadun, A family of M\"obius invariant
2-knot energies.- Geometric topology (Athens, GA, 1993), 
Amer. Math. Soc., Providence, 1997, 235--258.

\item{\ref{B}} William Browder, Diffeomorphisms of 1-connected 
manifolds. - 
Trans. Amer. Math. Soc. {\bf 128} (1967), 155--163. 

\item{\ref{BFHW}} Steve Bryson, Michael H.Freedman, Zheng-Xu He and
Zhenghan Wang, M\"obius invariance of knot energy.- Bull. Amer. Math.
Soc., {\bf 28}(1993), 99-103.

\item{\ref{BZ}} Gerhard Burde and Heiner Zieschang, Knots. - 
Walter de Gruyter \& Co., Berlin/New York, 1985.

\item{\ref{Gl}} Herman Gluck, The embedding of two-spheres in the
four-sphere. - Trans. Amer. Math. Soc. {\bf 104} (1962), 308-333.

\item{\ref{Go}} Camerun McA. Gordon, Knots in the 4-sphere.- Comment. Math.
Helv. {\bf 51} (1976), 585-596.

\item{\ref{GLT}} Mikhael Gromov, Blaine Lawson and William Thurston,
Hyperbolic 4-manifolds and conformally flat 3-manifolds.- Publ. Math. IHES,
{\bf 68}(1988), 27-45.

\item{\ref{Ka}}  Yoshi Kamishima,  Conformally  flat   manifolds  whose
development maps  are not surjective.  - Trans. Amer.  Math.
Soc. {\bf 294} (1986), 607-624.

\item{\ref{Ki}} Shin'ichi Kinoshita, Alexander polynomials as isotopy
invariants.- Osaka Math. J. {\bf 10} (1958), 263-271.

\item{\ref{Ko}} Christos Kourouniotis,
Deformations of hyperbolic structures.- Math. Proc. Cambr. Phil. Soc.
 {\bf 98} (1985), 247-261.

\item{\ref{Ku}}  Nicolaas Kuiper, Hyperbolic 4-manifolds and tesselations.-
Publ. Math. IHES, {\bf 68}(1988), 47-76.

\item{\ref{Kul}} Ravi Kulkarni, Infinite regular coverings. - 
Duke Math. J. {\bf 45} (1978), 781--796.

\item{\ref{KP}}  Ravi Kulkarni  and  Ulrich Pinkall,  Uniformization of
geometric   structures   with   applications   to  conformal
geometry.  - Differential  Geometry, Pe\~niscola,  1985 (Lect.
Notes Math. {\bf1209}), Springer, 1986, 190-210.

\item{\ref{LS}} L.K. Lashof and J.L. Shaneson, Classification of knots
in codimension two.-
Bull. Amer. Math. Soc. {\bf 75 } (1969), 171-175.

\item{\ref{L}} J. Levine, Polynomial invariants of knots of
codimension two. - Ann. of Math., {\bf 84} (1966), 537-554.

\item{\ref{Ma}} Gregory A. Margulis, The isometry of closed manifolds of
constant
negative curvature with the same fundamental group. - Soviet Math. Dokl.
{\bf 11}(1970), 722-723.

\item{\ref{Mas}} Bernard Maskit, Kleinian groups. - Springer-Verlag, 1987.

\item{\ref{Mo}} George D. Mostow, Quasiconformal mappings  in n-space and
the rigidity  of hyperbolic space forms.  - Publ. Math. IHES
{\bf 34} (1968), 53-104.

\item{\ref{MS}}  John Morgan  and  Peter Shalen,  Valuations,  trees and
degenerations of hyperbolic structures.  - Ann. of Math. {\bf 120}
(1984), 401-476.

\item{\ref{Mo}} John  Morgan,   
     Trees and hyperbolic geometry. - Proc.
Intern. Congress of Math. at Berkeley, {\bf 1}, 1987, 590-597.

\item{\ref{Ro}} Dale Rolfson, Knots and links.- Publish or Perish
Press, Berkeley, 1976.

\item{\ref{Ri}} Igor Rivin, On geometry of convex polyhedra in
hyperbolic 3-space.- Ph.D. Thesis, Princeton, 1986.

\item{\ref{Sp}} Edwin H. Spanier, Algebraic topology.-Springer-Verlag,
1966.

\item{\ref{Su}} Dennis Sullivan, Quasiconformal homeomorphisms and
dynamics, Part 2: Structural stability implies hyperbolicity for
Kleinian groups.-Acta Math. {\bf 155}(1985), 243-260.

\item{\ref{Suz}} Shin'ichi Suzuki, Knotting problems of 2-spheres in 
4-sphere.-  Math. Sem.
Notes Kobe Univ. {\bf 4} (1976), 241--333.

\item{\ref{Sw}} Gadde A.~Swarup, A note on higher dimensional knots.- 
Math. Zeit. {\bf 121} (1971), 141--144.

\item{\ref{T}} William Thurston, Hyperbolic  structures on 3-manifolds, I:
deformations of acylindrical manifolds. - Ann. of Math. {\bf 124}
(1986), 203-246.

\item{\ref{Tu}} Pekka Tukia, On isomorphisms of geometrically finite
Kleinian groups.- Publ. Math. IHES {\bf 61}(1985), 171-214.

\item{\ref{TV}} Pekka Tukia and Jussi V\"ais\"al\"a, Quasisymmetric
embeddings of metric spaces.- Ann. Acad. Sci. Fenn. Ser.A I Math.
{\bf 5}(1980), 97-114.

\enddocument